\documentclass[12pt]{amsart}
\usepackage{amsmath, amssymb, euscript, hyperref}

\newtheorem{thm}{Theorem}[section]
\newtheorem{lem}[thm]{Lemma}
\newtheorem{prop}[thm]{Proposition}
\theoremstyle{definition}
\newtheorem{defn}[thm]{Definition}

\newtheorem{rem}[thm]{Remark}
\newtheorem{cor}[thm]{Corollary}

\newcommand{\blackboard}[1]{\ensuremath{\mathbb{#1}}}

\newcommand{\N}{\blackboard{N}}

\newcommand{\Z}{\blackboard{Z}}

\newcommand{\R}{\blackboard{R}}

\newcommand{\F}{\ensuremath{\mathbf{F}}}

\begin{document}

\bibliographystyle{amsplain}

\address{Azer Akhmedov, Department of Mathematics,
North Dakota State University,
Fargo, ND, 58102, USA}
\email{azer.akhmedov@ndsu.edu}

\begin{center} {\bf On groups of homeomorphisms of the interval with finitely many fixed points } \end{center}

\medskip

 \begin{center} Azer Akhmedov \end{center}

\bigskip

{\bf Abstract:} {\Small We strengthen the results of \cite{A1}, consequently, we improve the claims of \cite{A2} obtaining the best possible results. Namely, we prove that if a subgroup $\Gamma $ of $\mathrm{Diff}_{+}(I)$ contains a free semigroup on two generators then $\Gamma $ is not $C_0$-discrete. Using this, we extend the H\"older's Theorem in $\mathrm{Diff}_{+}(I)$ classifying all subgroups where every non-identity element has at most $N$ fixed points. In addition, we obtain a non-discreteness result in a subclass of homeomorghisms which allows to extend the classification result to all subgroups of $\mathrm{Homeo}_{+}(I)$ where every non-identity element has at most $N$ fixed points.} 

\vspace{1cm}

\section{Introduction}

 Throughout this paper we will write $\mathbf{\Phi }$ (resp. $\mathbf{\Phi }^{\mathrm{diff}}$) to denote the class of subgroups of $\Gamma \leq \mathrm{Homeo}_{+}(I)$ (resp. $\Gamma \leq \mathrm{Diff}_{+}(I)$) such that every non-identity element of $\Gamma $ has finitely many fixed points. Let us point out immediately that any subgroup of $\mathrm{Diff}_{+}^{\omega }(I)$ - the group of orientation preserving analytic diffeomorphisms of $I$ - belongs to $\mathbf{\Phi }$. In fact, many of the major algebraic and dynamical properties of subgroups of $\mathrm{Diff}_{+}^{\omega }(I)$ is obtained solely based on this particular property of analytic diffeomorphisms having only finitely many fixed points. Interestingly, groups in $\mathbf{\Phi }$ may still have both algebraic and dynamical properties not shared by any subgroup of $\mathrm{Diff}_{+}^{\omega }(I)$. In particular, not every group in $\mathbf{\Phi }$ is conjugate to a subgroup of $\mathrm{Diff}_{+}^{\omega }(I)$.

 \medskip

  For a non-negative integer $N\geq 0$, we will also write $\mathbf{\Phi }_N$ (resp. $\mathbf{\Phi }^{\mathrm{diff}}_N$) to denote the class of subgroups of $\Gamma \leq \mathrm{Homeo}_{+}(I)$ (resp. $\Gamma \leq \mathrm{Diff}_{+}(I)$) such that every non-identity element of $\Gamma $ has at most $N$ fixed points in the interval $(0,1)$. For $f\in  \mathrm{Homeo}_{+}(I)$, $Fix(f)$ will denote the set of fixed of points of $f$ in $(0,1)$. 
  
  \medskip
  
  Characterizing $\mathbf{\Phi }_N$ for an arbitrary $N$ has been a major problem. In the case of $N=0$, H\"older's Theorem states that any subgroup of $\mathbf{\Phi }_0$ is Abelian, while in the case of $N=1$, Solodov's Theorem states \footnote{Solodov's result is unpublished but three independent proofs have been given by Barbot \cite{B}, Kovacevic \cite{K}, and Farb-Franks \cite{FF1}} that any subgroup of $\mathbf{\Phi }_1$ is metaabelian, in fact, it is isomorphic to a subgroup of $\mathrm{Aff}_{+}(\mathbb{R})$ - the group of orientation preserving affine homeomorphisms of $\mathbb{R}$.  
  
  \medskip   
  
  It has been proved in \cite{A1} that, for $N\geq 2$, any subgroup of $\mathbf{\Phi }^{\mathrm{diff}}_N$ of regularity $C^{1+\epsilon }$ is indeed solvable, moreover, in the regularity $C^2$ we can claim that it is metaabelian. The argument there fails short in complete characterization of subgroups of $\mathbf{\Phi }^{\mathrm{diff}}_N, N\geq 2$ even at these increased regularities. 
  
  \medskip
  
  In \cite{N1}, Navas gives a different proof of this result for groups of analytic diffeomorphisms, namely, it is shown that any group in $\mathbf{\Phi }^{\mathrm{diff}}_N$ of class $C^{\omega }$ is necessarily metaabelian. 
  
  \medskip
  
  In this paper, we provide a complete characterization of the class $\mathbf{\Phi }^{\mathrm{diff}}_N$ and $\mathbf{\Phi }_N$ for an arbitrary $N$. Our first main result is the following theorem.
  
  \medskip
  
 \begin{thm}\label{thm:main} Let $\Gamma \leq \mathrm{Diff}_{+}(I)$ be an irreducible subgroup, and $N\geq 0$ such that every non-identity element has at most $N$ fixed points. Then $\Gamma $ is isomorphic to a subgroup of $\mathrm{Aff}_{+}(\mathbb{R})$; moreover, if $\Gamma $ is not Abelian, then   $\Gamma $ is conjugate to a subgroup of $\mathrm{Aff}_{+}(\mathbb{R})$ 
 \end{thm} 
  
 \medskip
 
  In other words, any irreducible subgroup of $\mathbf{\Phi }^{\mathrm{diff}}_N$ is isomorphic to an affine group. Indeed, we show that, for $N\geq 2$, any irreducible subgroup of $\mathbf{\Phi }^{\mathrm{diff}}_N$ belongs to $\mathbf{\Phi }^{\mathrm{diff}}_1$! The conjugacy claim  is also proved in \cite{FF1} for non-Abelian subgroups $G\leq  \mathrm{Diff}_{+}^2(\R )$ where every non-identity element has at most one fixed point (see Theorem 1.5 there). In the proof, the authors establish that there is no wandering interval of the form $\theta ^{-1}(y) = [x_0, x_1], x_0 < x_1$ where $\theta = \nu (0, x)$ is the semi-conjugacy defined there. Our conjugacy claim (i.e. improving from semi-conugacy to a conjugacy)  also follows from the non-existence of wandering intervals, and in our case, this claim follows from the non-discreteness result (Theorem A$'$ of the next section). 

\medskip 

Let us point out that there exist meta-abelian examples  (communicated to the author by A.Navas; a certain non-standard representation of the Baumslag-Solitar group $BS(1,2)$ in $\mathrm{Homeo}_{+}(I)$) which shows that the class $\mathbf{\Phi }_N$ is indeed strictly larger than the class $\mathbf{\Phi }_1$, for $N\geq 2$. However, these groups are semi-conjugate to an affine group. Indeed, this semi-conjugacy has also been shown to hold  in the work of J.Carnevale \cite{C} when $N=2$; it is conjectured there that the semi-conjugacy holds for any $N$. In this paper, we discuss the cases of $N = 2$ and $N=3$ separately and provide proofs; then we prove the case of general $N$.  Let us remind that in \cite{A2} we also discuss the cases of $N = 1, 2, 3, 4$ separately and provide elementary proofs, but these proofs are based on non-discreteness results (hence, our results there are applicable to $C^{1+\epsilon }$-diffeomorhisms). For our proofs in the cases of $N\in \{2,3\}$ in this paper, we do not use non-discreteness and instead use a weak substitute for it. To discuss the cases of general $N$,  in the last section of the current paper, we obtain $C^0$-non-discreteness result for the groups in the class  $\mathbf{\Phi }_N$ which satisfy certain dynamical transitivity condition. Namely, we prove the following theorem. 
   
\medskip 

\begin{thm} \label{non-discretePhi} Let $\Gamma \in \mathbf{\Phi }_N$ for some $N\geq 1$ such that for any two open intervals $J_1, J_2$ with $\overline {J_1}\subset (0,1), \overline {J_2}\subset (0,1)$, there exists $g\in \Gamma $ such that $g(J_1)\subset J_2$. Then $\Gamma $ is non-discrete in the $C^0$ metric.

\end{thm}

 The proof of the main theorem  in \cite{A2} uses the bi-order $<$ (which can be defined also in the continuous category with all the used properties still holding) and the $C^0$-non-discreteness result. Thus we obtain the following result. 

 \begin{thm} \label{thm:classify} Let $\Gamma \in \mathbf{\Phi }_N$ for some $N\geq 1$. Then $\Gamma $ is semi-conjugate in $\mathrm{Homeo}_{+}(\mathbb{R})$  to a subgroup of $\mathrm{Aff}_{+}(\mathbb{R})$ 

\end{thm}

  \bigskip

   \section{$C_0$-discrete subgroups of $\mathrm{Diff}_{+}(I)$: strengthening the results of \cite{A1}} 
   
   \medskip
   
   The main results of \cite{A2} are obtained by using Theorems B-C from \cite{A1}. Theorem B (Theorem C) states that a non-solvable (non-metaabelian) subgroup of $\mathrm{Diff}_{+}^{1+\epsilon }(I)$ (of $\mathrm{Diff}_{+}^{2}(I)$) is non-discrete in $C_0$ metric. Existence of $C_0$-small elements in a group provides effective tools in tackling the problem. Theorems B-C are obtained by combining Theorem A in \cite{A1} by the results of Szekeres, Plante-Thurston and Navas. Theorem A states that for a subgroup $\Gamma \leq \mathrm{Diff}_{+}(I)$, if $[\Gamma , \Gamma ]$ contains a free semigroup in two generators then $\Gamma $ is not $C_0$-discrete. In the proof of Theorem A, the hypothesis that the generators of the free semigroup belong to the commutator subgroup $[\Gamma , \Gamma ]$ is used only to deduce that the derivatives of both of the generators at either of the end points of the interval $I$ equal 1. Thus we have indeed proved the following claim: {\em Let $\Gamma \leq \mathrm{Diff}_{+}(I)$ be a subgroup containing a free semigroup in two generators $f, g$ such that either $f'(0) = g'(0) = 1$ or $f'(1) = g'(1) = 1$. Then $\Gamma $ is not $C_0$-discrete, moreover, there exists non-identity elements in $[\Gamma , \Gamma ]$ arbitrarily close to the identity in $C_0$ metric}.    

  \medskip  

 In this section, we make a simple observation which strengthens Theorem A further, namely, the condition  ``$[\Gamma , \Gamma ]$ contains a free semigroup" can be replaced altogether with ``$\Gamma $ contains a free semigroup" (i.e. without demanding the extra condition ``either $f'(0) = g'(0) = 1$ or $f'(1) = g'(1) = 1$".
 
 \medskip
 
 \begin{thm}[Theorem A$'$] Let $\Gamma \leq \mathrm{Diff}_{+}(I)$ be a subgroup containing a free semigroup in two generators. Then $\Gamma $ is not $C_0$-discrete, moreover, there exists non-identity elements in $[\Gamma , \Gamma ]$ arbitrarily close to the identity in $C_0$ metric.   
  \end{thm}

  \medskip
  
  In the proof of Theorems B-C, if we use Theorem A$'$ instead of Theorem A we obtain the following stronger versions.
  
  \medskip
  
  \medskip
  
  \begin{thm} [Theorem B$'$] If a subgroup $\Gamma \leq \mathrm{Diff}_{+}^{1+\epsilon }(I)$ is $C_0$-discrete then it is virtually nilpotent.
  \end{thm}
  
  \medskip
  
  \begin{thm} [Theorem C$'$] If a subgroup $\Gamma \leq \mathrm{Diff}_{+}^2(I)$ is $C_0$-discrete then it is virtually Abelian.
  \end{thm}
  
  \medskip
  
   Theorem A$'$ is obtained from the proof of Theorem A by a very slight modification. Let us first assume that $\Gamma $ is irreducible, i.e. it has no fixed point on $(0,1)$. Let $f, g\in \Gamma $ generate a free semigroup on two generators. If $f'(0) = g'(0) = 1$ or $f'(1) = g'(1) = 1$ then the claim is already proved in \cite{A1}, otherwise, without loss of generality we may assume that $f'(1) < 1$ and $g'(1) < 1$. 
   
   \medskip
   
   Let also $\epsilon , N, \delta , M, \theta $ be as in the proof of Theorem A in \cite{A1}, except we demand that  $1 < \theta _N < \sqrt[8N]{1.9}$ (instead of $1 < \theta _N < \sqrt[2N]{2}$), and instead of the inequality $\frac{1}{\theta _N} < \phi '(x) < \theta _N$, we demand that $$\displaystyle \max _{x,y\in [1-\delta , 1]}(\frac{\phi '(x)}{\phi '(y)})^8 < \theta _N,$$ where $\phi \in \{f, g, f^{-1}, g^{-1}\}$. In addition, we also demand that for all $x\in [1-\delta , 1]$, we have $f(x) > x$ and $g(x) > x$. 
   
 \medskip
 
 Then we let $W = W(f,g), \alpha , \beta \in \Gamma $ be as in the proof of Theorem A. We may also assume that (by replacing $(\alpha , \beta )$ with $(\alpha \beta, \beta \alpha )$ if necessary), $\alpha '(0)  = \beta '(0) = \lambda < 1$. 
  
  \medskip
   
 Now, for every $n\in \mathbb{N}$, instead of the set  $$\mathbb{S}_n = \{U(\alpha ,\beta )\beta \alpha  \ | \ U(\alpha ,\beta ) \ \mathrm{is \ a \ positive \ word \ in} \ \alpha , \beta  \ \mathrm{of \ length \ at \ most} \ n\}$$
 
 we consider the set  $$\mathbb{S}_n' = \{U(\alpha ,\beta )\beta \alpha \in \mathbb{S}_n \ | \ \mathrm{sum \ of \ exponents \ of} \ \alpha \ \mathrm{in} \ U(\alpha ,\beta )  \ \mathrm{equals} \ [\frac{n}{2}]\}$$
 
  \medskip

  Previously, we had the crucial inequality $|\mathbb{S}_n| \geq 2^n$ for all $n$ but now we have the inequality $|\mathbb{S}_n'| \geq (1.9)^n$ for sufficiently big $n$. Let us also observe that, for any interval $J$ in $(1-\delta , 1)$, and for all $g\in \mathbb{S}_n'$, we will have the inequality $|g(J)| < \lambda ^n (\theta _N)^{\frac {1}{8}n}$. Then for some sufficiently big $n$ the following conditions hold: 
 
 \medskip
 
 (i) there exist $g_1, g_2\in \mathbb{S}_n$ such that $g_1\neq g_2$, and $$|g_1W(x_i)-g_2W(x_i)| < \frac{1}{\sqrt[2N]{1.9}^n} , 1\leq i\leq N-1,$$
 
 \medskip
 
 (ii) $M^{2m+4}(\theta _N)^{4n}\frac{1}{\sqrt[2N]{1.9}^n} < \epsilon $,
 
 \medskip
 
 where $x_i = \frac{i}{N}, 0\leq i\leq N$.
 
 \bigskip
 
 The rest of the proof goes exactly the same way by replacing $\mathbb{S}_n$ with $\mathbb{S}_n'$: letting again $h_1 = g_1W, h_2 = g_2W$, we obtain that $|h_1^{-1}h_2(x) - x| < 2\epsilon $ for all $x\in [0,1]$. Since $\epsilon $ is arbitrary, we obtain that $\Gamma $ is not $C_0$-discrete. On the other hand, by definition of $\mathbb{S}_n'$ we have $h_1^{-1}h_2\in [\Gamma , \Gamma ]$. If $\Gamma $ is not irreducible then it suffices to observe that there exists only finitely many intervals $I_1, \ldots , I_m$ in $(0,1)$ such that $\Gamma $ fixes the endpoints of $I_j$ but no other point inside $I_j$, moreover, $\displaystyle \sum _{1\leq j\leq m}|I_j| > 1-2\epsilon $  $\square $
   
   \bigskip
   
   \section{Extension of H\"older's Theorem in $\mathrm{Diff}_{+}(I)$}
   
   \medskip
   
   Let us point out that the following theorem follows from the proof of Theorem 0.1 and Theorem 0.2 in \cite{A2}. 
   
   \begin{thm}\label{thm:non-discrete} Let $N\geq 0$ and $\Gamma $ be an irreducible group in  $\mathbf{\Phi }^{\mathrm{diff}}_N$ such that $[\Gamma ,\Gamma ]$ contains diffeomorphisms arbitrarily close to the identity in $C_0$ metric. Then $\Gamma $ belongs to  $\mathbf{\Phi }^{\mathrm{diff}}_1$ thus it is isomorphic to a subgroup of the affine group  $\mathrm{Aff}_{+}(\mathbb{R})$.  
   \end{thm}

   \medskip
   
   The method of \cite{A2} does not allow to obtain a complete classification of subgroups of $\mathbf{\Phi }^{\mathrm{diff}}_N$ primarily because existence of non-discrete subgroups in $\mathrm{Diff}_{+}^{1+\epsilon }(I)$ (in $\mathrm{Diff}_{+}^{2}(I)$) is guaranteed only for non-solvable (non-metaabelian) groups. Within the class of solvable (metaabelian) groups the method is inapplicable.  
   
   \medskip

   Now, by Theorem A$'$, we can guarantee the existence of non-discreteness in the presence of a free semigroup. On the other hand, the property of containing a free semigroup on two generators is generic only in $C^{1+\epsilon }$ regularity; more precisely, any non-virtually nilpotent subgroup of $\mathrm{Diff}_{+}^{1+\epsilon }(I)$ contains a free semigroup on two generators. Just in $C^1$-regularity, $\mathrm{Diff}_{+}(I)$ has many non-virtually nilpotent subgroups (e.g. subgroups of intermediate growth) without free semigroups. (see \cite{N2})  
  
  \medskip
  
  The next proposition indicates a strong distinctive feature for groups of $\mathbf{\Phi }$, and supplies free semigroups for all non-Abelian subgroups in $\mathbf{\Phi }_N, N\geq 1$. 
  
  \medskip
  
  \begin{prop}\label{prop:semigroup} Any subgroup in class $\mathbf{\Phi }$ is either Abelian or contains a free semigroup on two generators.
  \end{prop}
  
  \medskip
  
  \begin{cor} For any $N\geq 0$, a subgroup of $\mathbf{\Phi }_N$ is either Abelian or contains a free semigroup.
  \end{cor}
  
  \medskip
  
  \begin{rem} Let us point out that any group $\Gamma $ in $\mathbf{\Phi }$ is bi-orderable. A bi-order can be given as follows: for $f, g\in \Gamma $, we let $f<g$ iff $f(x) < g(x)$ in some interval $(0, \delta )$. Proposition \ref{prop:semigroup} shows that the converse is far from being true, i.e. not every finitely generated bi-orderable group embeds in $\mathbf{\Phi }$. For example, it is well known that every finitely generated torsion-free nilpotent group is bi-orderable hence it embeds in $\mathrm{Homeo}_{+}(I)$ (by the result of \cite{FF2} it embeds into $\mathrm{Diff}_{+}(I)$ as well); on the other hand, a finitely generated nilpotent group does not contain a free semigroup on two generators.
  \end{rem}

 \medskip
   
  We need the following well known notion.
 
 \medskip
  
  \begin{defn} Let $f, g\in \mathrm{Homeo}_{+}(I)$. We say the pair $(f,g)$ is {\em crossed} if there exists a non-empty open interval $(a, b)\subset (0,1)$ such that one of the homeomorphisms fixes $a$ and $b$ but no other point in $(a,b)$ while the other homeomorphism maps either $a$ or $b$ into $(a,b)$.
  \end{defn}
  
 \medskip
 
  It is a well known folklore result that if $(f,g)$ is a crossed pair then the subgroup generated by $f$ and $g$ contains a free semigroup on two generators (see \cite{N3}).         
  
  \medskip
  
{\bf Proof of Proposition \ref{prop:semigroup}.} We may assume that $\Gamma $ is irreducible. If $\Gamma $ acts freely then by H\"older's Theorem it is Abelian and we are done. Otherwise, there exists a point $p\in (0,1)$ which is fixed by some non-identity element $f$ of $\Gamma $. Since $\Gamma $ is not irreducible, there exists $g$ which does not fix $p$. Let 
$p_{-}$ be the biggest fixed point of $g$ less than $p$, and $p_{-}$ be the smallest fixed point of $g$ bigger than $p$. If at least one of the points $p_{-}, p_{+}$ is not fixed by $f$ then either the pair $(f,g)$ or $(f^{-1},g)$ is crossed. 

\medskip

 Now assume that both  $p_{-}, p_{+}$ are fixed by $f$. Without loss of generality we may also assume that $g(x) > x$ for all $x\in (p_{-}, p_{+})$. Let $q_{-}$ be the smallest fixed point of $f$ bigger than $p_{-}$, and $q_{+}$ be the biggest fixed point of $f$ smaller than $p_{+}$. (we have $q_{-}\leq q_{+}$ but it is possible that $q_{-}$ equals $q_{+}$). Then there exists $n\geq 1$ such that $g^n(q_{-}) > q_{+}$. Then either the pair $(g^nfg^{-n}, f)$ or the pair  $(g^nf^{-1}g^{-n}, f)$ is crossed (in the interval $(a,b) = (q_{+}, p_{+})$). $\square $ 

  \medskip 

  We will need the following result which slightly strengthens Proposition \ref{prop:semigroup}.

  \medskip 

  \begin{prop} \label{prop:semigr2} A non-abelian subgroup in class $\mathbf{\Phi }$ contains two elements $h_1, h_2$ generating a free semigroup such that the smallest fixed point of $h_1$ in $(0,1)$ is smaller than the smallest fixed point of $h_2$ in $(0,1)$.   
 \end{prop}

  \begin{proof} Indeed, in the proof of Proposition \ref{prop:semigroup}, we can take $p$ to be the smallest fixed point of $f$ in $(0,1)$. Then if $p_{-} > 0$, we can take $h_1 = g, h_2 = f$, otherwise we let $h_1 = f, h_2 = g$.
  \end{proof}
 
 \vspace{1cm}
 
 \section{Semi-archimedian groups}
 
 It is a well known fact that any subgroup of $\mathrm{Homeo}_{+}(\mathbb{R})$ is left-orderable. Conversely, one can realize any countable left-orderable group as a subgroup of $\mathrm{Homeo}_{+}(\mathbb{R})$ (see \cite{N3}). Despite such an almost complete and extremely useful characterization of left-orderable groups, when presented algebraically (or otherwise) it can be difficult to decide if the group does admit a left order at all.
 
\medskip

  Let $G$ be a group with a left order $<$. $G$ is called {\em Archimedean} if for any two positive elements $f, g\in G$, there exists a natural number $n$ such that $g^n > f$. In other words, for any positive element $f$, the sequence $(f^n)_{n\geq 1}$ is {\em strictly increasing} and {\em unbounded}.\footnote{In a left-orderable group $G$, we say a sequence $(g_n)_{n\geq 1}$ is bounded if there exists an element $g$ such that $g^{-1} < g_n < g$ for all $n\geq 1$.}  It is a classical result, proved by H\"older, that Archimedean group are necessarily Abelian, moreover, they are always isomorphic to a subgroup of $\mathbb{R}$. In fact, the notion of Archimedean group arises very naturally in proving the fact that any freely acting subgroup of  $\mathrm{Homeo}_{+}(\mathbb{R})$ is Abelian, first, by showing that such a group must be Archimedean, and then, by a purely algebraic argument (due to H\"older), proving that {\em Archimedean $\Rightarrow $ Abelian.}       
   
  \medskip
  
  It turns out one can generalize the notion of Archimedean groups to obtain algebraic results of similar flavor for subgroups of $\mathrm{Homeo}_{+}(\mathbb{R})$ which do not necessarily act freely but every non-trivial element has at most $N$ fixed points. Let us first consider the following property.
  
  \begin{defn} Let $G$ be a group with a left order $<$. We say $G$ satisfies property $(P_1)$ if the following condition holds: for  a natural number $M$ and elements $g, \delta \in G$, if the sequence $(g^n)_{n\geq 1}$ is increasing but bounded, and $\delta g^k>g^m$ for all $k, m > M$, then for all $k\geq M$ either the sequence $(g^n \delta g^k)_{n\geq 1}$ or the sequence $(g^{-n} \delta g^k)_{n\geq 1}$ is increasing and unbounded.  
  \end{defn} 
  
  Every Archimedean group clearly satisfies property $(P_1)$ but there are non-archimedean groups too with property $(P_1)$. In fact, it is easy to verify that the metaabelian affine group $\mathrm{Aff}_{+}(\mathbb{R})$ with the following very natural order does satisfy property $(P_1)$ while not being Archimedean: for any two maps $f, g\in \mathrm{Aff}_{+}(\mathbb{R})$ we say $f<g$ iff either $f(0) < g(0)$ or $f(0) = g(0), f(1) < g(1)$. 
  
  \medskip
  
  An Archimedean group can be viewed as groups where powers of positive elements {\em reach infinity}. In groups with property $(P_1)$, the power of a positive element reaches infinity perhaps after an extra {\em arbitrarily small one time push}, namely if $g\in G$ is positive and $(g^n)_{n\geq 1}$ is still bounded, then for every $\delta $ where $\delta g^m > g^k$ for all sufficiently big $m, k$, either the sequence ${g^n\delta g^m}_{n\geq 1}$ or the sequence ${g^{-n}\delta g^m}_{n\geq 1}$ reaches the infinity. Thus groups with property $(P_1)$ can be viewed as generalization of Archimedean groups. We would like to introduce even a more general property $(P_N)$ for any $N\geq 1$. (Archimedean groups can be viewed as exactly the groups with property $(P_0)$).
  
  \medskip
  
  \begin{defn} \label{defn:baba} Let $G$ be a group with a left order $<$, and $N$ be a natural number. We say $G$ satisfies property $(P_N)$ if the following condition holds: for  a natural number $M$, elements $g, \delta _1, \dots , \delta _{N-1}\in G$, and the numbers $\epsilon _1, \dots , \epsilon _{N-1}\in \{-1,1\}$ if for all $i\in \{1, \dots , N-1\}$ and for all $k_1, \dots , k_{i-1}, k_i\geq M , \epsilon _1, \dots , \epsilon _i\in \{-1,1\}$, 
  
  (i) the sequence $(g^{\epsilon _in}\delta _{i-1}g^{\epsilon _{i-1}k_{i-1}}\dots \delta _{1}g^{\epsilon _{1}k_1})_{n\geq 1}$ is bounded from above,
  
  and 
  
  (ii) $\delta _i g^{\epsilon _ik_i}\delta _{i-1}g^{\epsilon _{i-1}k_{i-1}}\dots \delta _{1}g^{\epsilon _{1}k_1} > g^{\epsilon _ik_i}\delta _{i-1}g^{\epsilon _{i-1}k_{i-1}}\dots \delta _{1}g^{\epsilon _{1}k_1}$ 
  
  then, for some $\epsilon _N\in \{-1,1\}$, the sequence $(g^{\epsilon _Nn}\delta _{N-1}g^{\epsilon _{N-1}k_{N-1}}\dots \delta _{1}g^{\epsilon _{1}k_1})_{n\geq 1}$ is unbounded from above.    
  \end{defn}
  
  \medskip
 
  \begin{rem} Similarly, in groups with property $(P_N)$ the power of a positive element may not necessarily reach the infinity but does so after some $N$ arbitrarily small pushes (by $\delta _1, \dots , \delta _N$). Namely, one considers the sequences $$g^n, g^{\pm n}\delta _1g^n, g^{\pm n}\delta _2g^{\pm n}\delta _1g^n, \dots , g^{\pm n}\delta _N\dots g^{\pm n}\delta _1g^n$$ and one of them reaches infinity as $n\to \infty $. 
  \end{rem}
  
  \medskip
    
  \begin{rem} In the case of $N = 0$, the existence of elements $g_1, \delta _1, \dots , g_{N-1}, \delta _{N-1}$ is a void condition, and one can state condition $(P_0)$ as the existence of an element $g_0$ such that $g_0^n$ is unbounded; thus groups with property $(P_0)$ are exactly the Archimedean groups. 
  \end{rem}

  \begin{defn} A left ordered group $G$ is called {\em semi-Archimedean} if it satisfies property $(P_N)$ for some $N\geq 0$.  
  \end{defn}

  We will need the following result about semi-Archimedean groups:
  
  \begin{prop} \label{thm:archimed} Let $G$ be a countable semi-Archimedean group. Then $G$ has a realization as a subgroup of  $\mathrm{Homeo}_{+}(\mathbb{R})$ such that every non-identity element has at most $N$ fixed points.
  \end{prop}
  
   {\bf Proof.} For simplicity, we will first prove the proposition for $N = 1$.

   \medskip
   
  If there exists a smallest positive element in $\Gamma $ then, necessarily, $\Gamma $ is cyclic and the claim is obvious. Let $g_1, g_2, \dots $ be all elements of $\Gamma $ where $g_1 = 1$. We can embed $\Gamma $ in $\mathrm{Homeo}_{+}(\mathbb{R})$ such that the sequence $\{g_n(0)\}_{n\geq 1}$ satisfies the following condition: $g_1(0) = 0$,  and for all $n\geq 1$, 
  
  \medskip
  
  (i) if $g_{n+1} > g_i$ for all $1\leq i\leq n$, then $g_{n+1}(0) = \max \{g_i(0) \ | \ 1\leq i\leq n\} + 1$,
  
  \medskip
  
  (ii) if $g_{n+1} < g_i$ for all $1\leq i\leq n$, then $g_{n+1}(0) = \min \{g_i(0) \ | \ 1\leq i\leq n\} - 1$,
  
  \medskip
  
  (iii) if $g_i < g_{n+1} < g_j$ and none of the elements $g_1, \dots , g_n$ is strictly in between $g_i$ and $g_j$ then $g_{n+1}(0) = \frac{g_i(0)+g_j(0)}{2}$.
  
  \medskip 
  
  Then, since there is no smallest positive element in $\Gamma $, we obtain that the orbit $O = \{g_n(0)\}_{n\geq 1}$ is dense in $\mathbb{R}$. This also implies that the group $\Gamma $ for any point $p\in O$ and for any open non-empty interval $I$, there exists $\gamma \in \Gamma $ such that $\gamma (p)\in I$.  
  
  \medskip
  
  Now assume that some element $g$ of $\Gamma $ has at least two fixed points. Then for some $p, q$ we have $Fix(g)\cap [p,q] = \{p,q\}$. Without loss of generality, we may also assume that $p > 0$ and $g(x) > x$ for all $x\in (p,q)$. By density of the orbit $\{g_n(0)\}_{n\geq 1}$, there exists $f\in \Gamma $ such that $f(0) \in (p,q)$. Then, for sufficiently big $n$, we have $\delta = g^{-n}f$ has a fixed point $r\in (p,q)$, moreover, $\delta (x) > x$ for all $x\in (p,r)$. 
  
  \medskip
  
  Then $g^{\epsilon n}$ does not reach infinity for any $\epsilon \in \{-1,1\}$, in fact, $g^{\epsilon n}(0) < p$ for all $n\geq 1, \epsilon \in \{-1,1\}$. Then $\{g^{\epsilon _1n}\delta g^{\epsilon k}\}_{n\geq 1}$ does not reach infinity for any $k\in \Z, \epsilon , \epsilon _1\in \{-1,1\}$. Contradiction.  
   
   \medskip
   
  To treat the case of general $N\geq 1$, let us assume that some element $g\in \Gamma $ has at least $N+1$ fixed points. Then there exists open intervals $I_1 = (a_1, b_1), \dots , I_{N+1} = (a_{N+1}, b_{N+1}$ such that $a_1 < b_1 \leq a_2 < b_2 \leq  \dots \leq a_N < b_N \leq a_{N+1} < b_{N+1}$ and $\{a_1, b_1, \dots , a_{N+1}, b_{N+1}\subset Fix(g)$. By density of the orbit $O$, there exist elements $\delta _1, \dots , \delta _N$ such that  $\delta _i(b_i)\in I_{i+1}, 1\leq i\leq N$. Then for the appropriate choices of $\epsilon _1, \dots , \epsilon _{N-1}\in \{-1,1\}$ and for sufficiently big $k_1, \dots , k_{N_1}$, conditions (i) and (ii) of Definition \ref{defn:baba} hold, while for any $\epsilon _N\in \{-1,1\}$, the sequence $(g^{\epsilon _Nn}\delta _{N-1}g^{\epsilon _{N-1}k_{N-1}}\dots \delta _{1}g^{\epsilon _{1}k_1})_{n\geq 1}$ is bounded from above because it lies in $I_{N+1}$. $\square $
  
  \medskip
  
  \begin{rem} Let us emphasize that in this section we did not make an assumption that the groups belong to the class $\mathbf{\Phi }$.
  \end{rem} 
   
 \vspace{1cm}

\bigskip 

 \section {Classification of subgroups of $\mathbf{\Phi }_N$ up to a semi-conjugacy: cases of small $N$}
 
 This section is devoted to Theorem \ref{thm:classify}. In the language of semi-Archimedean groups, this theorem implies the following proposition. 

  \begin{prop} Any semi-Archimedean group is semi-conjugate in $\mathrm{Homeo}_+(I)$ to an affine group.  
\end{prop} 
 
 \medskip 
 
We recall some ingredients that we have used in \cite{A2}. Let $\Gamma $ be in lass $\mathbf{\Phi }$. We will write $$F(\Gamma) = \displaystyle \mathop{\sup }_{g\in \Gamma \backslash \{1\}} |Fix(g)|,$$ thus $\F(\Gamma )\in \N \cup \{\infty \}$.  A fixed point $p\in (0,1)$ of $\gamma \in \Gamma $ is called {\em tangential}, if there exists $\delta > 0$ such that either $\gamma (x) \geq x$ for all $x\in (p-\delta , p+\delta )$ or $\gamma (x) \leq x$ for all $x\in (p-\delta , p+\delta )$; otherwise $p$ is called {\em non-tangential}. 

 For $f, g\in \Gamma $, we say 

 (i) $g$ {\em separates} $f$ if $Fix(g)\neq \emptyset , |Fix(f)| \geq 2$ and the set $Fix(g)$ lies in the interval $(a,b)$ where $a, b$ are the smallest and the second smallest fixed points of $f$ respectively; 

 (ii) $g$ {\em precedes} $f$, if $Fix(f)\neq \emptyset , Fix(g)\neq \emptyset  $ and for all $x\in Fix(g), y\in Fix(f)$ we have $x < y$.
 
 \medskip 
 
 We will assume that $N\geq 2$ as the case of $N = 0$ and $N=1$ follow from the theorems of H\"older and Solodov respectively. 
 
 \medskip 
 
 In the class $\mathbf{\Phi }$ we can introduce a bi-order which we have used in \cite{A2}: For $f, g \in \Gamma $ , we let $g < f$ if $g(x) < f(x)$ near zero. 
 
 \medskip 
 
  For a general bi-ordered group $\Gamma $ (not necessarily an embedded subgroup of $\mathrm{Homeo}_+(\R )$) with (an abstract) bi-order $<$ we will need the following notions: For $f\in \Gamma $ we let $|f| = f$ is $f$ is positive and $|f| = f^{-1}$ otherwise. If $f, g \in \Gamma $, we say $g$ is {\em infinitesimal w.r.t.} $f$ (or $g$ is infinitesimal w.r.t. $f$ near zero) if $|g|^n < |f|$ for every $n\in \mathbb{Z}$; in this case we write $g << f$. If neither of $f, g$ is infinitesimal w.r.t. the other, we say $f$ and $g$ are {\em comparable}. If an element $\gamma \in \Gamma $ is not infinitesimal w.r.t. any other element, then we say $\gamma $ is a dominant element of $\Gamma $. Notice that if $\Gamma $ is finitely generated, then one of the generators is necessarily a dominant element. Also, a conjugate of a dominant element is also dominant, and all dominant elements lie outside of the commutator subgroup (since a commutator $[x,y]$ is infinitesimal w.r.t. both $x$ and $y$)   

  \medskip 
  
  For $f\in \Gamma $, we also write $\Gamma _f = {\gamma  \in \Gamma  : \gamma << f}$. Notice that if $\Gamma $ is finitely generated with a fixed finite symmetric generating set, and $f$ is the biggest generator, then $\Gamma _f$ is a normal subgroup of $\Gamma $ (which we call {\em the infinitesimal subgroup of $\Gamma $}), moreover, $\Gamma /\Gamma _f$ is Archimedean, hence Abelian. Consequently, we have $[\Gamma , \Gamma] \leq  \Gamma _f$.  

 \medskip 

  Now we return back to a subgroup $\Gamma $ from the class $\mathbf{\Phi }$. Similar to the bi-order $<$, we need a bi-order $<'$ by letting $f<'g$ if $f(x)<g(x)$ bear 1. If $f$ is infinitesimal w.r.t. $g$ w.r.t. both $<$ and $<'$ then we say $f$ is infinitesimal w.r.t. $g$ near 0 and near 1.

  \medskip

We need to recall some basic notions and facts about the subgroups of  $\mathrm{Homeo}_+(I)$. A subgroup $\Gamma $ is called {\em irreducible} if  $\Gamma $ has no global fixed point in $(0,1)$. Borrowing from \cite{A2}, we call a subgroup $\Gamma $ {\em dynamically transitive} if for all $p\in (0,1)$ and for any interval $I\subseteq (0,1)$, there exists $\gamma \in \Gamma $ such that $\gamma (p)\in I$.   

\medskip 

It is useful to view subgroups of  $\mathrm{Homeo}_+(I)$ also as subgroups of $\mathrm{Homeo}_+(\mathbb{R})$. Let $G\leq \mathrm{Homeo}_+(\mathbb{R})$. A non-empty closed subset $M\subseteq \mathbb{R}$ is called a minimal set of $G$ if it is $G$-invariant and does not contain any proper closed $G$-invariant subset. The following result is folklore:

 \begin{thm} \label{thm:minimalset} Let $G\leq \mathrm{Homeo}_+(\mathbb{R})$. Then either there is no minimal set or for any minimal set $M$ one of the following mutually exclusive cases hold: 
 
  (1) $M$ is discrete; in this case, $M\subseteq \{x: \forall g \ \mathrm{with} \ Fix(g)\neq \emptyset , x\in Fix(g)\}$ and the set $\{x: \forall g \ \mathrm{with} \ Fix(g)\neq \emptyset , x\in Fix(g)\}$ is a union of minimal sets.
  
  (2) $M$ is a perfect nowhere dense set; in this case, $M$ is a unique minimal set and $M\subseteq O(x) := \{g(x) : g\in G\}$ for every $x\in \mathbb{R}$
  
  (3) $M = \mathbb{R}$.
  
  In addition, if $G$ is finitely generated then it admits a minimal set. 
  \end{thm}

\medskip 

  In cases (2) and (3) of the theorem when the minimal set is unique, we will write $M = \mathcal{M}(G)$. 
  
  \medskip 
  
  In case (2), the group will be semi-conjugate to a group $G_1\leq \mathrm{Homeo}_+(\mathbb{R})$ with $\mathcal{M}(G) = \mathbb{R}$. When we view a subgroup of  $\mathrm{Homeo}_+(\mathbb{R})$  as a subgroup of $\mathrm{Homeo}_+(I)$, then instead of  $\mathcal{M}(G) = \mathbb{R}$ we write $\mathcal{M}(G) = (0,1)$. 

  \medskip 
  
  Let us also remind another folklore result about the classification of subgroups of $\mathrm{Homeo}_+(\mathbb{R})$. 

   \begin{thm} \label{thm:proximal} Let $G\leq \mathrm{Homeo}_+(\mathbb{R})$ without a global fixed point . Then one of the following cases holds: 
   
   (1) $G$ admits an invariant Radon measure on $\mathbb{R}$, in which case the group surjects onto $\mathbb{Z}$ - the additive group of integers.
   
(2) $G$ does not preserve a Radon measure on $\mathbb{R}$, but it is semi-conjugate to a minimal action that commutes with integer translations. 

 (3)  for all bounded, open intervals $I$ and $J$, there exists $g\in G$ such that $g(I) \subset J$.

   \end{thm}

  Let us point out that an irreducible subgroup from the class $\mathbf{\Phi }_N$, viewed as a subgroup of $\mathrm{Homeo}_+(\mathbb{R})$, never preserves a Radon measure and never commutes with the translations unless it acts freely, i.e. it belongs to $\mathbf{\Phi }_0$. In the latter case, by H\"older's Theorem, the group is Abelian. 

  \medskip 
  
   From the descriptions above we immediately obtain the following
  
  \begin{thm} \label{thm:transitive} $G\leq \mathrm{Homeo}_+(I)$ be a non-Abelian irreducible subgroup from the class $\mathbf{\Phi }$. Then $G$ is semi-conjugate to a subgroup $G_1\leq \mathrm{Homeo}_+(I)$  satisfying the property that for any two open intervals $J_1, J_2$ with $\overline {J_1}\subset (0,1), \overline {J_2}\subset (0,1)$, there exists $g\in G_1$ such that $g(J_1)\subset J_2$. In particular, $\mathcal{M}(G_1) = (0,1)$
  \end{thm}
  
  Theorem \ref{thm:transitive} motivates the following notion
  
  \begin{defn} A subgroup $G\leq \mathrm{Homeo}_+(I)$ is called {\em proximal} if for any two open intervals $J_1, J_2$ with $\overline {J_1}\subset (0,1), \overline {J_2}\subset (0,1)$, there exists $g\in G_1$ such that $g(J_1)\subset J_2$. 
  
  \end{defn} 

   Let us note that proximal subgroups are necessarily dynamically transitive, and dynamically transitive subgroups are irreducible.

   \medskip 
   
   As a corollary of Theorem \ref{thm:transitive}, we obtain that 
  
  \begin{prop} \label{thm:transitive2} Let $\Gamma \in \mathbf{\Phi }_N$ for some $N\geq 1$. Then, $\Gamma $ is either Abelian or  semi-conjugate to a proximal subgroup $\Gamma _1\in \mathbf{\Phi }_N$.
  \end{prop}

 We will need the following result

  \begin{prop} \label{thm:arrange} Let $\Gamma \leq \mathrm{Homeo}_+(I)$ be an irreducible non-solvable group with $\mathcal{M}(\Gamma ) = (0,1)$, $F(\Gamma ) = N\geq 1$ and $0 < a < b < 1$. Then there exists $f, h\in \Gamma $ with $|Fix(h)|\geq 2 $ such that $f$ and $h$ agree near 0 and near 1, $h << f$ near 0 and near 1, $f$ and $h$ generate a free semigroup, $f$ has at least two non-tangential fixed points, and for all $x\in Fix(f), y\in Fix(h), y < a < b < x$. Moreover, if $\Gamma $ is finitely generated, then $f$ can be chosen as a dominant element near 0 and near 1.
  \end{prop}
  
  \medskip
  
  \begin{proof} 
  
   For $N=1$ the claim follows from Solodov's Theorem. Let $N\geq 2$. 

   \medskip 
   
   Let us make an important observation that if $\Gamma \in \Phi _N$ is irreducible, then so is the commutator subgroup $[\Gamma , \Gamma ]$ (and hence any subgroup  $\Gamma ^{(n)}$ in the derived series of $\Gamma $). On the other hand, if $\Gamma  \in \Phi _N$ is irreducible and (as assumed) non-solvable, then it has a finitely generated subgroup which is also irreducible and non-solvable. 

   \medskip 
   
   Thus we can assume that $\Gamma $ is finitely generated and non-solvable. Let $S$ be a finite generating set of $\Gamma $. There exists $f_{left}\in \Gamma \backslash [\Gamma , \Gamma ]$ such that $Fix(f_{left})\neq \emptyset $ and $f_{left}$ is dominant near 0. Indeed, otherwise we have $f_1\in \Gamma \backslash [\Gamma , \Gamma ], f_2\in [\Gamma , \Gamma ], f_1\in S$ such that $f_i(x) > x, i\in \{1,2\}$ near zero,  $Fix(f_1) = \emptyset , Fix(f_2)\neq \emptyset $, $g<<f_1$ near zero for all $g\in [\Gamma , \Gamma ]$ (in particular, $f_2<<f_1$ near zero) and $f_1$ is a dominant element. Then for sufficiently big $n$, we can take $f_{left} := f_2^{-n}f_1$; then  $f_{left}$ has a fixed point, (for all $n$) $f_{left}\in \Gamma \backslash [\Gamma , \Gamma ]$ and $f_{left}$ is a dominant element near zero. Similarly, we can claim that  there exists $f_{right}\in \Gamma \backslash [\Gamma , \Gamma ]$ such that $Fix(f_{right})\neq \emptyset $ and $f_{right}$ is dominant near 1. Then necessarily there exists $\bar{f}\in \{f_{left}^nf_{right}^m : m, n\in \Z\}$ such that $Fix(\bar{f})\neq \emptyset $, $\bar{f}$ is positive and it is dominant both near zero and near one. 

\medskip 

   By H\"older's Theorem, we also have $h_1\in [\Gamma , \Gamma ]$ such that $h_1(x) > x$ near zero and $|Fix(h_1)|\geq 2 $. 
   
\medskip 

   Since $\Gamma $ is irreducible and $\Gamma \in F(N)$, $[\Gamma , \Gamma ]$ is also irreducible. Then there exists $g_1, g_2\in [\Gamma , \Gamma ]$ such that for all $y\in Fix(g_1h_1g_1^{-1}), x\in Fix(g_2\bar{f}g_2^{-1})$ we have $y < a < b < x$. On the other hand, $g_1h_1g_1^{-1} << \bar{f}$ and $\bar{f}$ is comparable with $g_2\bar{f}g_2^{-1}$ hence $g_1h_1g_1^{-1}<<g_2\bar{f}g_2^{-1}$. Then for sufficiently big $m$, the homeomorphisms $g_1h_1^mg_1^{-1}$ and $f: = g_2\bar{f}g_2^{-1}$ are crossing (hence their big enough powers  generate a free semigroup). Thus we can take $h = g_1h_1^mg_1^{-1}$ where $m$ is sufficiently big. 

   \medskip 

   $f$ and $h$ will agree near 0, indeed they will be both positive. If they do not agree near 1, then we can replace the pair $(\Gamma , \Gamma ^{(1)})$ with $(\Gamma ^{(1)}, \Gamma ^{(2)})$ and choose $f_{new}\in \Gamma ^{(1)}, h_{new}\in \Gamma ^{(2)}$ such that $f_{new}<<f, h_{new}<<f_{new}$ satisfying the required conditions, namely, $|Fix(h_{new})|\geq 2 $, $f_{new}$ and $h_{new}$ are both positive (hence agree near 0),  $h_{new} << f_{new}$ near 0, $f_{new}$ has at least two non-tangential fixed points, and for all $x\in Fix(f_{new}), y\in Fix(h_{new}), y < a < b < x$. If $f_{new}$ and $h_{new}$ do not agree near 1, then the pair $f$ and $h_{new}$ will agree both near 0 and near 1. Thus we can take $f$ and a conjugate $\psi h_{new}\psi ^{-1}$ to satisfy the required conditions. 
   
  \end{proof}

\medskip

  Proposition \ref{thm:arrange} allows to classify subgroups of $\Phi _N$ for $N = 2$ and $N=3$. 

  \begin{prop} \label{N=2} Let $\Gamma \leq \mathrm{Homeo}_+(I)$  with $\mathcal{M}(\Gamma ) = (0,1)$ and $\Gamma \in \Phi _2$. Then $\Gamma $ is semi-conjugate to an affine group.    
  \end{prop}

  \begin{proof} It suffices to show that $F(\Gamma ) \leq 1$ (then the claim of our proposition follows from Solodov's Theorem).  If $\Gamma $ is  solvable, then the claim follows from Plante's result \cite{P} that any solvable subgroup  $G\leq \mathrm{Homeo}_+(\R)$ where $Fix(\gamma )$ is discrete is semi-conjugate to an affine action. If $\Gamma $ is not solvable, then by Proposition \ref{thm:arrange}, we have positive $f, h\in \Gamma $ such that $h$ has at least two fixed points,  $f, h$ agree near zero and near one, $h<<f$ near zero and near one, and $h$ precedes $f$ (i.e. for all $x\in Fix(f), y\in Fix(h), y < x$). Let $a$ be the smallest fixed point of $h$. Then for sufficiently big $n$, $h^{-n}f$ has at least two non-tangential fixed points in $(0,a)$ and at least one non-tangential fixed point in $(a,1)$. Thus $|Fix(h^{-n}f)| \geq 3$. Contradiction.

  \end{proof}

   \medskip 

   The proof of the above proposition crucially uses Proposition \ref{thm:arrange} which, in its turn, makes use of the observations that for any $N\geq 1$, if $\Gamma \in \Phi _N$ is irreducible and non-solvable, then so is the $n$-th commutator subgroup $\Gamma ^{(n)}$ for any $n\geq 1$. On the other hand, if $\Gamma  \in \Phi _N$ is irreducible and (as assumed) non-solvable, then it has a finitely generated subgroup which is also irreducible and non-solvable. In the proof of the following proposition, we also need to make use of the proximality of $\Gamma $, and although it is true that proximality of $\Gamma \in \Phi _N$ implies proximality of $\Gamma ^{(n)}$ for any $n\geq 1$, but we have a technical difficulty in allowing us the assumption that $\Gamma $ is also finitely generated.

   \medskip

 \begin{prop} \label{N=3} Let $\Gamma \leq \mathrm{Homeo}_+(I)$  with $\mathcal{M}(\Gamma ) = (0,1)$ and $\Gamma \in \Phi _3$. Then $\Gamma $ is semi-conjugate to an affine group.    
  \end{prop}

   \begin{proof} Again, we will show that $F(\Gamma ) \leq 1$. If $\Gamma $ is solvable, the claim follows from Plante's result \cite{P} so we will assume that $\Gamma $ is not solvable. 

    \medskip 

    By Proposition \ref{thm:arrange} and by the proximality of $\Gamma $ (hence of $\Gamma ^{(n)}$ for any $n\geq 0$), there exist positive $h_1, \psi _1\in \Gamma ^{(3)}, f_1\in \Gamma ^{(2)}\backslash \Gamma ^{(3)}$ with $|Fix(h_1)|\geq 2 $ such that $f_1$ and $h_1$ agree near 0 and near 1,  $h_1 << f_1$ near 0 and near 1, $f_1$ has at least two non-tangential fixed points, $h_1$ precedes $f_1$, and $\psi _1$ separates $h_1$.  

    \medskip 

     Let $G_1$ be the subgroup of $\Gamma ^{(1)}$ generated by a finite subset $\Omega _1$ where the element $f_1$ can be written as a product of commutators of elements of $\Omega _1$ and the elements $h_1, \psi _1$ can be written as product of double commutators of elements of $\Omega _1$. 

     \medskip 

     Applying Proposition \ref{thm:arrange} and using the proximality of $\Gamma $ again, we can choose positive $h_2, \psi _2\in \Gamma ^{(2)}, f_2\in \Gamma ^{(1)}\backslash \Gamma ^{(2)}$ with $|Fix(h_2)|\geq 2 $ such that $f_2$ is comparable with all elements of $\Omega _1$, $f_2$ and $h_2$ agree near 0 and near 1,  $h_2 << f_2$ near 0 and near 1, $f_2$ has at least two non-tangential fixed points, $h_2$ precedes $h_1$ and $f_1$ precedes $f_2$ (hence $h_2$ precedes $f_2$), and $\psi _2$ separates $h_2$.

     \medskip 

     We now let $G_2$ be the subgroup of $\Gamma $ generated by a finite subset $\Omega _2$ where the element $f_2$ can be written as a product of commutators of elements of $\Omega _2$ and the elements $h_2, \psi _2$ can be written as product of double commutators of elements of $\Omega _2$. Then, applying Proposition \ref{thm:arrange} and using the proximality of $\Gamma $ again, we can choose positive $f_3\in \Gamma \backslash \Gamma ^{(1)}$  such that $f_3$ is comparable with all elements of $\Omega _1\cup \Omega _2$, $f_3$ has at least two non-tangential fixed points and $f_2$ precedes $f_3$. 

      \medskip 

      Then, necessarily, we have $f\in \{f_2, f_3\}, \psi \in \{\psi _1, \psi _2\}, h\in \{h_1, h_2\}$ such that $\psi $ separates $h$, $h$ precedes $f$, $h<<f$ and $\psi <<f$ near 0 and near 1, and $f$ and $h$ agree near 0 and near 1. We have $\psi (x) < x$ near zero. If $\psi (x) < x$ near 1, then we let $\phi _n := \psi ^{-n}h\psi ^{n}$; otherwise, we let  $\phi _n := \psi ^{n}h\psi ^{-n}$. Then, for sufficiently big $n$, $\phi _n$ has fixed points $a, b$ such that $Fix(f)\subset (a,b)$. Then for sufficiently big $k$, $\phi _n^{-k}f$ has at least four fixed points. Contradiction.
   \end{proof}
   
\bigskip 

 \section{Non-discreteness of expansive subgroups}

 We now start proving a non-discreteness result for proximal non-solvable subgroups of $\Phi _N, N\geq 1$. Let us remind that in the proof of main results of \cite{A2}, besides discreteness, all the arguments work in the continuous category as well. Thus to prove Theorem \ref{thm:classify}, we need need to show that non-solvable subgroups of  $\mathbf{\Phi }_N$ are semi-conjugate to a $C^0$-non-discrete subgroup. \footnote{Let us recall that a subgroup $G\leq \mathrm{Homeo}_{+}(I)$ is called $C^0$-discrete if there eists $\epsilon > 0$ such that $||g||_0 > \epsilon $ for all $g\in G\backslash \{1\}$ where $||g||_0 = \displaystyle \mathop {\max }_{0\leq x\leq 1}|g(x)-x|$.}

\medskip 

 In dealing with non-smooth homeomorphisms, one of the issues causing difficulty for us is that for a homeomorphism $h$, in the given interval $J$, we may not be able to find infinitely many conjugates of $h$ by a fixed diffeomorphism such that each conjugate maps some subinterval  $J'\subseteq J$ to an interval of length at least $c|J'|$ for a uniformly fixed constant $c > 0$. It turns out this "bad case" can be used in the direction of showing non-discreteness by the method already exploited in Section 2 and in \cite{A1}. 

 \medskip 

 \begin{prop} \label{prop:expanding} Let $\Gamma $ be a proximal $C_0$-discrete subgroup in class $\mathbf{\Phi }$, $f, h\in \Gamma $ generate a free semigroup with $f\in \mathrm{Diff}_{+}(I)$ and $d = \min Fix(h) > 0$ such that $1 < f'(x) < 1.1, f''(x) > 0$ and $x < h^4(x) < f(x)$ for all $x\in (0,d)$. Then for every non-empty interval $J = (a, b)\subset (0,d)$ with $f(J)\cap J = \emptyset $, there exists an infinite subset $A\subseteq \mathbb{N}$ such that for all $n\in A$ there exists a sub-interval $J_n\subseteq (f^{-3}(a),b)$ such that  $\displaystyle \mathop{\inf }_{n\in A}|J_n| > 0$ and $\displaystyle \mathop{\inf }_{n\in A}\frac{|f^{n}hf^{-n}(J_n)|}{|J_n|} > 0$.  
     
 \end{prop}

 \begin{proof} Let $\epsilon > 0$ be such that $||\gamma || > \epsilon $ for all $\Gamma \backslash \{1\}$. Let $N$ be a natural number such that $1/N < \epsilon $ and $x_i = i/N, 1\leq i\leq N-1$. Let also $n_0$ be a positive integer. 

  \medskip 
  
  Since $\Gamma $ is proximal, there exists $\gamma \in \Gamma $ such that $\gamma (x_i)\in J, 1\leq i\leq N-1$. Let $\delta = \displaystyle \mathop{\min }_{1\leq i\leq N-2}|\gamma (x_{i+1})-\gamma (x_i)|$. For all $n\geq 3, r\geq 1$, let $$S_n^r = \{W(f,h)\in P(f^{-1},h^{-1}) : |W| = n, W = W_1f^r, f^r\in Suffix(W)\}$$ where $P(f^{-1},h^{-1})$ is the set of positive words (viewed also as a subset of $\Gamma $) in the alphabet $\{f^{-1}, h^{-1}\}$ and $Suffix(W)$ is the set of all suffixes of $W$.

\medskip 

   Since $\Gamma $ is $C_0$-discrete, there exists $k\in \{1,\dots , n-1\}$ and an infinite subset $A\subseteq \mathbb{N}$ such that for all $n\in A$ and $r > n_0$ there exists $U_n, V_n\in S_n^{r}$ such that $$0 < V_n(x_k) - U_n(x_k) < \frac{1}{|S_n^r|} < \frac{1}{1.2^n} \ \mathrm{and}  \ |U_n^{-1}V_n(x_k) - x_k| > \delta .$$ By Mean Value Theorem (applied to iterations of $f$) and by inequality $1< f'(x) < 1.1$ for all $x\in (0,d)$, there exists a suffix $hV_n'$ of $V_n$ of length at least $n_0$ such that $|h(V_n'(I_n))| \geq |V_n'(I_n)|$ where $I_n = (U_n(x_k), V_n(x_k))$. Let us also assume that the length $|hV_n'|$ is minimal. Then by the inequality $|U_n^{-1}V_n(x_k) - x_k| > \delta $, for some $m\geq 1$, $f^m(h(V_n'(I_n)))\subset (f^{-3}a,b)$ and $|f^m(h(V_n'(I_n)))| > \frac{\delta }{2}$. Then, again by Mean Value Theorem $|f^mhf^{-m}(K_{m,n})| > \frac{1}{2}|K_{m,n}| > \frac{\delta }{4}$ where $K_{m,n} = f^m((V_n'(I_n)))$. \footnote{Let us observe that if $I', I'' = h(I')$ are subintervals in some interval $(p,q)$ with $f^2(p) > q, f^m(q) < d$ and $|I''| \geq |I'|$, then by Mean Value Theorem $\frac{|f^m(I'')|}{|f^m(I')|} = \displaystyle \mathop{\Pi }_{j=1}^{m}\frac{f'(f^{j-1}x)}{f'(f^{j-1}y)} \frac{|I''|}{|I'|}$ for some $x\in I', y\in J''$. Since $f'$ is monotone in $(0, \frac{1}{2})$, we obtain that $\displaystyle \mathop{\Pi }_{j=1}^{m}\frac{f'(f^{j-1}x)}{f'(f^{j-1}y)} \frac{|I''|}{|I'|} \leq (\displaystyle \mathop{\max } _{s,t\in (0,d)}\frac{f'(t)}{f'(s)})^2\frac{|I''|}{|I'|} < (1.1)^2\frac{|I''|}{|I'|} < (1.3)\frac{|I''|}{|I'|}$. Now it remains to let $I' = V_n'(I_n)$.}
   
   \medskip 
   
   We let $n_1 = m$ and $J_{n_1} = K_{m,n}$; so $|J_{n_1}| > \frac{\delta }{2}$ and $\frac{|f^{n_1}hf^{-n_1}(J_{n_1})|}{|J_{n_1}|} > \frac{\delta }{4}$. By choosing $n$ bigger $n_1$, this time we choose $n_2$ and $J_{n_2}$ and continue the process; in choosing $n_{i+1}$ we choose a suffix $hV_n'$ of length at least $n_i$. Then we let $A = \{n_1, n_2, \dots \}$.  

 \end{proof} 
 
 \medskip 

   The above proposition motivates the following notion

   \begin{defn} A family $\mathcal{F}$ of orientation preserving homeomorphisms of $I=[0,1]$ is called {\em uniformly semi-expansive on the interval $J\subseteq (0,1)$} if there exists a constant $c > 0$ such that for all $f\in \mathcal{F}$, there exists a sub-interval $J'\subseteq J$  such that $|J'| > c$ and $|f(J')| > c|J'|$. If $J = (0,1)$ then we simply say that $\mathcal{F}$ is uniformly semi-expansive.        
   \end{defn}

\medskip 

  Let us observe that for $f\in \mathrm{Diff}_{+}(I)$, the singleton family $\mathcal{F} = \{f\}$ is uniformly semi-expansive where one can take $c = \displaystyle \mathop{\min }_{x\in [0,1]} |f'(x)|$.. If $f\in \mathrm{Diff}_{+}^2(I)$ and $f''(x) < 0$ for all $x\in (0,1)$, then the family  $\mathcal{F} = \{f^nhf^{-n} : n\in \N \}$ is uniformly semi-expansive provided  $h$ is also smooth and  $x \leq h(x) \leq f(x)$ for sufficiently small $x$.  
  
 \medskip 

  Proposition \ref{prop:expanding} indeed establishes that the family $\{h_n : n\in A\}$ is uniformly semi-expansive in the interval $(0, 1/2)$. However, this result is still not sufficient for us to achieve a desired non-discreteness result. In the next proposition we will assume conditions which are more than what is provided by Proposition \ref{prop:expanding} (i.e. we have not established all the assumptions of the next proposition). To forumate these condtions, we will need the following modified notion

  \medskip 

   \begin{defn} Let $A\subseteq \N$ be an infinite subset. A family $\{\phi _n\}_{n\in A}$ in $\mathrm{Homeo}_{+}(I)$ is called {\em weakly uniformly expansive} on the interval $J\subseteq (0,1)$ if there exists a constant $c > 0$ such that for any sub-interval $J_1\subseteq J$ there exists a sub-interval $J_2\subseteq J_1$ such that $|\phi _n(J_2)| > c|J_2|$ for all but finitely many $n\in A$.  
        
    \end{defn}
    
    Let us note that the notion of {\em weakly uniformly expansive} is neither stronger nor weaker than the notion of {\em uniformly semi-expansive}. It would also interesting to compare weakly uniformly expansiveness with the notion of expansive set from the next section.

\medskip 

 \begin{prop} \label{thm:nondiscrete} Let $\Gamma $ be a proximal group belonging to the class $\mathbf{\Phi }$, $f, h\in \Gamma $ generate a free semigroup with $f\in \mathrm{Diff}_+(I)$, $d = \min Fix(h) > 0$ such that $x < h^4(x) < f(x)$ for all $x\in (0,d)$ and $h$ is uniformly semi-expansive in $(0,d)$. Then $\Gamma $ is non-discrete in the $C_0$ metric.
  \end{prop}

  \begin{proof}  

 Notice that $\Gamma $ is necessarily irreducible. Let $\epsilon > 0$, $m$ be a positive integer such that $\frac{1}{m} < \epsilon /4$ and $x_i = \frac{i}{m}, 1\leq i\leq m-1$. 

\medskip

  By conjugating $\Gamma $ in $\mathrm{Homeo}_+(I)$, we may assume that 

  \medskip 
  
  (c1) $f$ is smooth of class $C^2$, i.e. $f, f^{-1}\in C^2[0,1]$;

 \medskip 

  (c2) $1/2$ is the smallest fixed point of $f$ and $d$ is the smallest fixed point of $h$ where $0 < d < 1/4$;

  \medskip 

  (c3) $f''(x) < 0$ for all $x\in (0,1/2)$;

  \medskip 

  (c4) $1.1 > f'(x) > 1$ for all $x\in (0,1/4)$ and $1 > f'(x) > 0.9$ for all $x\in (1/4, 1/2)$; 

  \medskip 

 Let us notice that, as a result of conditions (c1)-(c4), by Mean Value Theorem, for all intervals $J_i\subset (0,d), 1\leq i\leq 2$, if $|J_1| < |J_2|$, then $|f(J_1)| < \frac{1.1}{0.9}|f(J_2)|$. 

\medskip 
  
  Since $\Gamma $ is proximal, there exist $g\in \Gamma $ such that $g(x_i) \in (0, d)$ for all $1\leq i\leq m-1$. Let $y_i = g(x_i), 1\leq i\leq m-1$ and $I_i = (y_i, y_{i+1}), 1\leq i\leq m-2$. Then, since $h$ is uniformly semi-expansive for all $k\geq 1$ and closed intervals $J_i\subset (0,1/2), 1\leq i\leq k$, there exists an infinite family $\{f^nhf^{-n} : n\in A\}$ which is uniformly semi-expansive in the intervals $J_i$. (Notice that by Proposition \ref{prop:expanding}, we can only claim that there exists an infinite subset $A\subseteq \mathbb{N}$ such that for all $n\in A$ and $k\in \{1, \dots , m-2\}$ there exists a sub-interval $J_n^{(k)}\subseteq (f^{-3}(y_k),y_{k+1})$ [instead of $J_n^{(k)}\subseteq (y_k,y_{k+1})$!] such that  $$\displaystyle \mathop{\inf }_{n\in A}|J_n^{(k)}| > 0 \ \mathrm{and} \ \displaystyle \mathop{\inf }_{n\in A}\frac{|f^{n}hf^{-n}(J_n^{(k)})|}{|J_n^{(k)}|} > 0.)$$  By taking $J_i = I_i, 1\leq i\leq m-2$, since $A$ is infinite, we obtain that  there exist $p, q\geq 1$ such that $h_p(y_i) \in (h_q(y_{i-1}), h_q(y_{i+1})$ for all $2\leq i\leq m-2$. Then $g^{-1}(h_p^{-1}h_q)g(x_i)\in (x_{i-1}, x_{i+1}), 2\leq i\leq m-2$. Thus $||g^{-1}(h_p^{-1}h_q)g||_0 < \epsilon $. This shows that $\Gamma $ is not discrete. Contradiction.

  \end{proof}

  \medskip 

  \begin{rem}
       By Proposition \ref{prop:semigr2}, there exist $f, h\in \Gamma $ generating a free semigroup such that $x < h(x) < f(x)$ near zero  and $\min Fix(h) < \min Fix(f)$. We can also satisfy all the assumptions of Proposition \ref{thm:nondiscrete} except we cannot satisfy the conditions about smoothness of $f$ and semi-expansiveness of $h$ simultaneously (by conjugating inside $\mathrm{Homeo}_{=}(I)$ it is straightforward to satisfy each of these conditions separately but it is a major problem to have them both arranged).
  \end{rem}

  \bigskip

  \section{Classification of subgroups of $\mathbf{\Phi }_N$: the case of general $N$}

\medskip 

 In this section we will outline another approach to prove Theorem \ref{non-discretePhi}. Since homoemorphisms are monotone functions, they are differentiable almost everywhere. This fact is encouraging in dealing with homemorphism groups of the interval (since a generic continuous function is nowhere differentiable, though it is true that functions  locally of bounded variation are almost everywhere differentiable as well), however, the set of points where differentiability holds may have a degenerate image. The famous Cantor stair function $\kappa :[0,1]\to [0,1]$ is differentiable in a set $S =[0,1]\backslash C$ of full Lebesgue measure where $C$ is the Cantor set, but the image $\kappa (S)$ will have zero measure. Such a wild behaviour may occur not just for increasing continuous functions, but also for homeomorphisms.  

\medskip 

  \begin{defn} An increasing continuous function $\phi :[a,b]\to [a,b]$ with $\phi (a) = a, \phi (b) = b$ is called {\em singular} if $\phi '(t) = 0$ for almost all $t\in [0,1]$. $h$ is called {\em regular} if $\phi '(t) > 0$ for almost all $t\in [0,1]$
  \end{defn}

\medskip 

   Cantor stair function is singular though it is not a homeomorphism; Salem \cite{S} has constructed the first example of a singular homeomorphism of the interval. Let us observe that if $\phi :[a,b]\to [a,b]$ is a singular homeomorphism, then for any diffeomorphism $f:[a,b]\to [a,b]$, the maps $f\phi , \phi f, f\phi f^{-1}$ are singular.  

\medskip 

   Given a subgroup $\Gamma \leq \mathrm{Homeo}_+(I)$ generated by two elements, we can conjugate the group to make one of the generators smooth. It is possible that the other generator is then conjugated to a singular homeomorphism. This possibility adds to the difficulty of the problem, although we want to emphasize that significant effort below is made even in the case when no element of $\Gamma $ is singular. 

   \medskip 

   Every homeomorphism $\phi \in \mathrm{Homeo}_{+}(I)$ is differentiable almost everywhere in $[0,1]$. Then for all $\phi \in \Gamma $ and $c > 0$, we have a measurable partition $[0,1] = A(\phi ,c)\sqcup B(\phi ,c)$ such that for all measurable $A\subseteq A(\phi ,c)$ we have $\lambda (\phi (A)) \leq c\lambda (A)$ and for all measurable $B\subseteq B(\phi ,c)$ we have $\lambda (\phi (B)) > c\lambda (B)$; here, $\lambda $ denotes the standard Lebesgue measure. Notice that for regular homeomorphisms, this decomposition not only exists but is also unique up to a measure zero set difference, moreover the set $B(\phi ,c)$ will have positive measure. We will call $B(\phi ,c)$ {\em an expansive set} of $\phi $ with constant $c>0$.

 \medskip 
 
  Let us also note that if $\phi $ is singular, then we also have a decomposition $[0,1] = C(\phi )\sqcup D(\phi )$ where $C(\phi ) = \{t\in [0,1] : f'(t) \ \mathrm{exists} \}$ (thus $C(\phi ) = \{t\in [0,1] : f'(t) = 0 \}$) and $D(\phi ) = [0,1]\backslash C(\phi )$. Then $C(\phi )$ has full measure and $D(\phi )$ is dense; moreover, $\phi (C(\phi ))\subseteq D(\phi ^{-1})$.

  \medskip 
  
  To proceed further, we need the following well known fact from basic functional analysis (for convenience of the reader, we provide its proof here). The motivation for this lemma comes from the fact that in a given interval, it is possible to have sets $E_n, n\geq 1$ and $c > 0$ such that $\lambda (E_n)\geq c$ for all $n\geq 1$ (here, $\lambda $ denotes the usual Lebesgue measure on $\mathbb{R} $) nevertheless, for any subsequence $E_{n_k}, k\geq 1$, the intersection $\displaystyle \mathop{\cap }_{k\geq 1}E_{n_k}$ has zero Lebesgue measure. Let us note that  if $E_n$ are intervals, then we would be able to claim that there exists a  subsequence $E_{n_k}, k\geq 1$ such that  the intersection $\displaystyle \mathop{\cap }_{k\geq 1}E_{n_k}$  has positive Lebesgue measure; we will use this fact in the sequel, but after discussing an application of the next 

\medskip 

  \begin{lem} \label{thm:weakconvergence} Let $\psi _n:[a,b]\to \R , n\geq 1$ be a sequence of functions in $L^2([a,b]$ weakly converging to $\psi \in L^2[a,b]$. Then the sequence of functions  $\xi _n(t) := \int _0^t\psi _n(t)dt, n\geq 1$ uniformly converges to $\xi (t) := \int _0^t\psi (t)dt$. 
  \end{lem}

  \begin{proof} First, from the equality $$|\xi _n(x)-\xi(x)| = |\int _a^x(\psi _n(x)-\psi (x))dx|= |\langle \mathbf{1}, \psi _n \rangle - \langle \mathbf{1}, \psi  \rangle|$$ we obtain that $\xi_n$ converges to $\xi $ pointwise. On the other hand, for all $a\leq x\leq y\leq b$, by H\"older's Inequality, we have $$|\xi_n(x)-\xi _n(y)| = |\int _x^y\psi_n(t)dt| \leq \sup _n||\psi _n||_{L^2}\sqrt{y-x}$$ where $||\psi _n||_{L^2}$ denotes the norm in $L^2[a,b]$. Then it follows from Arzela-Ascoli Theorem that $(\psi _n)$ has a uniformly converging subsequence. But since $(\psi _n)$ converges to $\psi $ pointwise, we obtain that the entire sequence $(\psi _n)$ (not just its subsequence) converges to $\psi $ uniformly.  
  \end{proof}

\medskip

  Now we present another set of arguments to prove the non-discreteness claim (as in Proposition \ref{thm:nondiscrete}) replacing the condition on uniform expansiveness with another condition on certain regularity in the behaviour of $h$.  Assume that $\Gamma $ is not solvable.  Let $\epsilon , m, f, h$ be as in the proof  Proposition \ref{thm:nondiscrete}; we also assume the conditions (c1)-(c4) from the previous section.  

\medskip

 We will make an assumption that the map $h$ is regular \ (1).

 \medskip 
 
 Let us observe that $B(f, 1/2) \supseteq (0,1/2)$. On the other hand, for sufficiently small $c>0$, $\lambda (B(h,c)\cap(0,d)) > 0$, hence  $\lambda (B(h,c)\cap(0,1/2)) > 0$. Let us assume that we have a uniform lower bound for $\lambda (B(h_n,c)\cap(0,1/2))$ along a subsequence, i.e. for some infinite subset $A = \{k_1, k_2, \dots \}\subset \N $ and $\delta > 0$, we have $$\displaystyle \mathop {\inf}_{n}\lambda (B(h_n,c)\cap(0,1/2)) > \delta  \ \ (2) $$ where $h_n = f^{k_n}hf^{-k_n}, n\geq 1$. 
 
 \medskip 

  Let $B_n = B(h_n,a)\cap (0,\frac{1}{2}), n\geq 1$. Then we have $\lambda (B_n)\geq \delta $ for all $n\geq 1$. The functions $\psi _n = \mathbf{1}_{B_n}$ belong to the unit ball of $L^2[0,1]$ hence by Banach-Alaoglu Theorem, the sequence $(\psi _n)$ has a weakly converging subsequence; we still denote with $(\psi _n)$. Then by Lemma \ref{thm:weakconvergence}, the sequence $\xi _n(x) := \int _0^x\psi _n(t)dt$ converges uniformly to $\xi (x) := \int _0^x\psi (t)dt$

\medskip 

 Now we claim that there exists a measurable $\Delta \subseteq [0,1]$ of positive measure that for any interval $J$ with $\lambda (J\cap \Delta) > 0$, $\xi $ is not constant on $J$. Indeed, otherwise there exists a disjoint union $\displaystyle \mathop{\sqcup }_{n\geq 1}J_n$ of countably many intervals such that $\lambda ([0,1]\backslash \displaystyle \mathop{\sqcup }_{n\geq 1}J_n) = 0$ and $\xi \equiv const$ on $J_n$ for all $n\geq 1$.

  \medskip

  Then there exists an integer $M\geq 1$ such that $\displaystyle \mathop{\sum }_{n\geq M}\lambda (J_n) < \frac{\delta }{4}$. On the other hand, since $(\xi _n)$ converges to $\xi $ uniformly, for all $1\leq i\leq M$, by taking $\epsilon = \frac{\delta }{2^{i+2}}$, we can claim that there exists a positive integer $N_i$ such that for all $n > N_i \lambda (B_n\cap J_i) < \frac{\delta }{2^{i+2}}$. Let $N = \max\{N_1, \dots , N_M\}$. Then for all $n > N$ we obtain $$\lambda (B_n) \leq \lambda (B_n\cap (\displaystyle \mathop{\cup }_{n\geq M}J_n)) + \displaystyle \mathop{\sum }_{1\leq k\leq M}\lambda (J_k\cap B_n) < \frac{\delta }{4} + \displaystyle \mathop{\sum }_{i=1}^M\frac{\delta }{2^{i+2}} < \delta $$ which is a contradiction.  

  \medskip

  Thus we established the existence of a measurable $\Delta \subseteq [0,1]$ with $\lambda (\Delta ) > 0$ such that on any interval  having intersection with $\Delta $ of positive measure $\xi $ is non-constant. Then for every measurable $D\subseteq \Delta $ with $\lambda (D) > 0$ we have $\displaystyle \mathop{\liminf }_{n\to \infty }\lambda (B_n\cap D) > 0 \ (\ast )$. Indeed, otherwise, by considering a linear continuous functional $\langle \mathbf{1}_D, \cdot \rangle $ on $L^2[0,1]$ and evaluating it on the sequence $\mathbf{1}_{B_n}$ (which we assumed converges weakly in $L^2[0,1]$) we obtain that $\displaystyle \mathop{\lim }_{n\to \infty }\lambda (B_n\cap D) = 0$. Then $\xi $ is constant on $D$. Then, since $\lambda (D) > 0$, $\xi $ is constant on some interval of positive length which is again a  contradiction.  

\medskip 

  Since $\Gamma $ is proximal and $\lambda (\Delta ) > 0$, there exist $y_i, 1\leq i\leq m-1$ and $g\in \Gamma $ such that for all $1\leq i\leq m-2$,  $y_i < x_i < y_{i+1}, g(y_i) \in \Delta $ and $\lambda ((g(y_i). g(y_{i+1}))\cap \Delta ) > 0$. Let $z_i = g(y_i), 1\leq i\leq m-1$. 

  \medskip 

  By $(\ast )$, there exists $\omega  > 0$ such that for all $n\geq 1$ and $1\leq i\leq n-1$, we have $|h_n(z_{i+1})-h_n(z_i)| \geq \omega \ (\ast \ast )$. Then, since inequality $(\ast \ast )$ holds for infinitely many $n$, there exist $p, q\geq 1$ such that $h_p(z_i) \in (h_q(z_{i-1}), h_q(z_{i+1})$ for all $2\leq i\leq m-2$. Then $g^{-1}(h_p^{-1}h_q)g(y_i)\in (y_{i-1}, y_{i+1}), 2\leq i\leq m-2$. Thus $||g^{-1}(h_p^{-1}h_q)g||_0 < \epsilon $. This shows that $\Gamma $ is not discrete which is a contradiction.

 \medskip 

    Now we need to discuss the general case, i.e. without the assumptions (1) and (2).

    \medskip 
    
    {\bf Proof of Theorem \ref{non-discretePhi}.} Since $\Gamma $ is proximal, it is not Abelian. Let us assume  $\Gamma $ is not discrete. Then there exits $\epsilon > 0$ such that $||g||_0 := \displaystyle \mathop {\max }_{x\in [0,1]}|g(x)-x| > \epsilon $ for all $g\in \Gamma \backslash \{1\}$. By Proposition \ref{thm:arrange}, we may assume that there exists positive $f, h\in \Gamma $ with $|Fix(h)|\geq 2 $ such that $f$ and $h$ agree near 0 and near 1, $h << f$ near 0 and near 1, $f$ and $h$ generate a free semigroup, $f$ has at least two non-tangential fixed points, and $h$ precedes $f$. Then $f^{-1}hf$ and $h$ generate an irreducible subgroup, hence again by Proposition \ref{thm:arrange}, for some reduced words $U, V\in \mathbb{F}_2$, the maps $U(f^{-1}hf, h)$ and $V(f^{-1}hf, h)$ generate a free semigroup. 

    \medskip 

     By conjugating $\Gamma $ in $\mathrm{Homeo}_+(I)$, we may assume that 

  \medskip 
  
  (d1) $f$ is smooth of class $C^2$, i.e. $f, f^{-1}\in C^2[0,1]$;

 \medskip 

  (d2) $x < h(x) < f(x)$ for all $x\in (0,1/2)$;

  \medskip 

  (d3) $f(x) = 1.1x$ for all $x\in (0,1/2)$.

  \medskip 

  Let $m \geq 1$ such that $\frac{1}{m} < \frac{\epsilon }{2}$ and $x_i = \frac{i}{m}, 1\leq i\leq m-1$ and $\gamma \in \Gamma $ such that $\gamma (x_{m-1}) < 1/2$ and $f(\gamma x_1) > \gamma x_{m-1}$ (such an element $\gamma $ exists because of proximality of $\Gamma $). Let also, for all $n, r\geq 0$, \\ $$\delta _n = \displaystyle \mathop{\min }_{1\leq i\leq m-2}|f^n(\gamma (x_{i+1}))-f^n(\gamma (x_{i}))|$$  and $$S_n^r = \{W(U_0,V_0)f^{-r}\gamma : W(U_0,V_0)\in P(U_0,V_0) : |W| = n\}$$ \\
  
  where $U_0 = U(f^{-1}h^{-1}f, h^{-1}), V_0 = V(f^{-1}h^{-1}f, h^{-1})$ and  $P(U_0,V_0)$ is the set of positive words (viewed also as a subset of $\Gamma $) in the alphabet $\{U_0, V_0\}$.

 \medskip 

   Then by conditions (d1)-(d3), we have $\delta _r = \frac{\delta _0}{1.1^r}$ for all $r\geq 0$. Since $h<<f$, for any fixed $n$, if $r$ is sufficiently big with respect to $n$, we have $\Omega _n^p \cap \Omega _n^q = \emptyset $ for all $p, q\geq r, |p- q| > 1$ where $$\Omega _j^i = \{g(x) : x\in (x_1, x_{m-1}), g\in S_j^i\} \ \mathrm{for \ all} \  i,j\geq 0.$$  On the other hand, since $|S_n^r| \geq 1.9^n$ for all $n, r\geq 0$ and  $\Gamma $ is not discrete, there exists $U, V\in S_n^r$ and $i\in \{1, \dots , m-1\}$ such that $$0 < U(x_l) - V(x_l) < (1.5)^{-n}, 1\leq l\leq m-1 \ \mathrm{and} \ |U^{-1}V(x_i)-x_i| > \frac{1}{m}.$$ Then there exists a constant $C > 0$ and $N\geq 1$ such that for all $n\geq N$, if $r$ is sufficiently big, then there exists an interval $J_r\subset \Omega _n^r, r\geq 1$ with $|J_r| > C(1.1)^{-r}$ and $|h(J_r)| > |J_r|$. Then  $|f^r(J_r)| > C$ and $|f^rh(J_r)| > C$ hence $|f^rhf^{-r}(f^r(J_r))| > C$. Let us note that, since $h<<f$, if $r$ is sufficiently big with respect to $n$, then  the interval $f^r(J_r)$ lies in the interval $(\min f^{-1}(\Omega _0^0), \max f(\Omega _0^0))$ where  $\Omega _0^0 = (\gamma (x_1), \gamma (x_{m-1}))$. 

   \medskip 

   Now, for all $k\geq 1, 1\leq i\leq m-1, 1\leq j \leq k$, let $$y_{k, i, j} = \gamma (x_i) + j\frac{\gamma (x_{i+1}) -\gamma (x_i)}{k}, 1\leq j\leq k.$$ Let $k\geq 1$ be fixed.  Then for all sufficiently big $n, r$, there exits $i\in \{1, \dots , m-1\}$ such that for all $1\leq j \leq k, 1\leq l\leq m-1$, there exists $U, V \in S_n^r$ such that $$0 < U(y_{k,l,j}) - V(y_{k,l,j}) < (1.5)^{-n}, 1\leq l\leq m-1$$ and $$|U^{-1}V(y_{k,i,j})-y_{k,i,j}| > \frac{1}{m}.$$ Then, by strengthening the result of the previous paragraph, we can claim that there exists a constant $C > 0$ (independent of $k$; indeed it suffices to take $C = \frac{\delta _0}{2m}$) and $N\geq 1$ such that for all  $n\geq N$ if $r$ is sufficiently big, then there exists an interval $J_r\subset \Omega _n^r$ with a partition $J_r = \displaystyle \mathop {\sqcup }_{1\leq j\leq k}J_{r,k,j}$ such that $$|J_r| > C(1.1)^{-r}, |J_{r,k,j}| = \frac{1}{k}|J_r| \  \mathrm{and} \ |h(J_{r,k,j})| > |J_{r,k,j}|.$$ Then $|f^rh(J_{r,k,j})| > \frac{C}{k}$ hence $|f^rhf^{-r}(f^r(J_{r,k,j}))| > \frac{C}{k}$. 
   
   \medskip 
   
   We also have $$|f^r(J_r)| > C, r\geq 1 \ \mathrm{and} \ |f^r(J_{r,k,j})| > \frac{C}{k}, r, k\geq 1, 1\leq j\leq k.$$ Then there exist an infinite subset $A\subseteq \N $ such that $\displaystyle \mathop {\cap }_{r\in A}f^r(J_r)$ contains an interval $J$ of positive length and the family $\{h_r : = f^rhf^{-r} : r\in A\}$ is weakly uniformly expansive on $J$. 

   \medskip 

    Since $\Gamma $ is proximal, we have $g\in \Gamma $ such that $g(\gamma (x_i))\in J$ for all $1\leq i\leq m-1$. Then, since the set $A$ is infinite, and the family $\{h_r : r\in A\}$ is weakly uniformly expansive on $J$, there exists $p, q\in A$ such that $||g^{-1}(h_p^{-1}h_q)g||_0 < \epsilon $. Contradiction. \ $\square $

    \bigskip 
    
    Comparing the above result with the non-discreteness result of \cite{A1}, we see that the assumption about smoothness of the group is dropped, but instead, our group is proximal and belongs to class $\mathbf{\Phi }_N$ for some $N\geq 1$.

 \end{document}